%% file: agt-3-46.tex
\documentclass{gtart}
\input agtout

\lognumber{46}
\volumenumber{3}
\volumeyear{2003}
\papernumber{46}
\published{31 December 2003}
\pagenumbers{1277}{1289}
\received{7 April 2003}
\revised{1 December 2003}
\accepted{9 December 2003}

\usepackage{amsmath,amssymb,amscd}
\let\Bbb\mathbf

\newtheorem{theorem}{Theorem}[section]

\newtheorem{proposition}[theorem]{Proposition}

\newcommand{\edim}{\hbox{\rm e-dim}}
\newcommand{\z}{{\Bbb Z}}
\newcommand{\q}{{\Bbb Q}}

\newcommand{\s}{{\Bbb S}}
\newcommand {\h}{{\check H}}

 \newcommand{\invlim}{\raisebox{-1ex}{$\stackrel{\hbox{lim}}{\leftarrow}$}}

\newcommand{\lo}{\longrightarrow}
\newcommand{\sm}{\setminus}
\newcommand{\tor}{{\rm Tor}}
\begin{document}

\title{Cell-like resolutions preserving cohomological\\dimensions}

 \author{Michael  Levin}
 \address{Department of Mathematics\\Ben Gurion University of the 
Negev\\P.O.B. 653\\Be'er Sheva 84105, ISRAEL}
  \email{mlevine@math.bgu.ac.il}

\begin{abstract}
We prove that for every compactum $X$   with
$\dim_\z X   \leq n  \geq 2$
there is a  cell-like  resolution  $r: Z \lo X$ from
 a compactum $Z$
  onto $X$  such that $\dim Z  \leq n$ and
   for   every integer $k$ and   every abelian group $G$
  such that  $\dim_G X \leq k \geq 2$   we have $\dim_G Z \leq k$.
 The latter property  implies that for every simply connected
 CW-complex $K$ such that    $\edim  X \leq K$  we also have
 $\edim Z \leq K$.
\end{abstract}
\asciiabstract{%
We prove that for every compactum X with dim_Z X <= n >= 2 there is a
cell-like resolution r: Z --> X from a compactum Z onto X such that
dim Z <= n and for every integer k and every abelian group G such that
dim_G X <= k >= 2 we have dim_G Z <=k.  The latter property implies
that for every simply connected CW-complex K such that e-dim X <=
K we also have e-dim Z <= K.}
\maketitle

\begin{section}{Introduction}
A space $X$  is always assumed to be separable metrizable. The cohomological dimension $\dim_G X$
of  $X$  with respect to an abelian group $G$ is   the least number $n$ such that $\h^{n+1}
(X,A;G)=0$ for every closed subset $A$ of $X$.
 It  was known long ago that
 $\dim X  =\dim_\z X$ if $X$ is finite dimensional.
 The first example of an infinite dimensional compactum (=compact metric space)
 with
 finite integral cohomological dimension  was
  constructed by Dranishnikov \cite{dr0}  in 1987.  In 1978 Edwards \cite{ed1, w1}
 discovered his celebrated resolution theorem revealing a  close  relation between
 $\dim_\z$ and $\dim$.  The    Edwards  resolution
 theorem says that a compactum of $\dim_\z \leq n$ can be obtained
 as the  image of  a cell-like map defined  on a compactum of  $\dim \leq n$.
   A  compactum    $X$ is cell-like if any map $f : X \lo K$ from $X$ to a CW-complex  $K$  is
   null-homotopic.
   A  map is cell-like if its fibers are cell-like. The reduced $\rm {\check C}$ech
 cohomology groups of  a cell-like compactum   are trivial with respect to any group $G$.

    The Edwards resolution theorem addresses only the case of integral cohomological dimension.
   It seems natural to  look for  generalizations of this theorem taking into consideration
   other abelian groups.     Indeed,  such an investigation has been of considerable interest
   in cohomological dimension theory.
   It mainly  went  in two directions.

   The first one is  to adjust   resolutions for a given
   group  $G$  replacing   cell-like maps by $G$-acyclic maps.     A map is $G$-acyclic if
    the reduced $\rm {\check C}$ech   cohomology groups modulo $G$ of the fibers are trivial.
     By the Vietoris-Begle  theorem a $G$-acyclic map
cannot raise the cohomological dimension $\dim_G$.
   Let us give two  examples of results of this type.

  \begin{theorem} {\rm \cite{dr1}}\qua
\label{dr1} Let $p$ be  a prime number and
 let $X$ be  a compactum with $\dim_{\z_p} X \leq n$.  Then there are a compactum
 $Z $  with
  $\dim Z \leq n$ and a $\z_p$-acyclic map $r : Z \lo X$  from $Z$ onto $X$.
\end{theorem}

    \begin{theorem}  {\rm \cite{l1}}\qua
\label{l1}
 Let $G$ be an abelian group and
 let $X$ be  a compactum with $\dim_G X \leq n$, $ n \geq 2$.  Then there are a compactum
 $Z $  with
 $\dim_G Z \leq n$ and  $\dim Z \leq n+1$ and a $G$-acyclic map $r : Z \lo X$  from $Z$ onto $X$.
\end{theorem}

      The other direction of investigation  is to construct
   cell-like resolutions     preserving cohomological dimensions with respect to several
   abelian groups.   Below are  some   results   of this type.

    \begin{theorem}  {\rm \cite{dr1.5}}\qua
    \label{dr1.5}
   Let $p$ be a prime  number and   let a   compactum $X$  be such that
   $\dim_{\z_p} X \leq n$ and $\dim_{\z[1/p] } X \leq n$, $n \geq 2$. Then there  are
  a compactum  $Z $  with
 $\dim Z \leq n+1$,     $\dim_{\z_p} Z \leq n$ and $\dim_{\z[1/p] } Z \leq n$
   and a cell-like  map $r : Z \lo X$  from $Z$ onto $X$.
 \end{theorem}
      \begin{theorem}  {\rm \cite{dr2}}\qua
\label{dr2}  Let   $\cal L$ be a subset  of the set of primes  and let
  $X$ be a compactum  such that  $\dim_\z X \leq n$ and  $ \dim_{\z_p} X \leq k$, $n < 2k-1$
 for every  $p \in \cal L$.  Then there are a compactum
 $Z $  with
 $\dim Z \leq n$ and  $\dim_{\z_p} Z \leq  k $
 for every $p \in \cal L$ and a cell-like  map $r : Z \lo X$  from $Z$ onto $X$.
\end{theorem}

\begin{theorem}     {\rm \cite{ko1}}\qua
\label{ko1}
Let $p,q$ be distinct prime numbers and let $n$ be an integer $>1$.
Then  for a compactum $X$ with $\dim_{\z_p} X \leq n$,
$\dim_{\z_{(q)}} X \leq n$ and $\dim_\z X \leq n+1$
there exist an $(n+1)$-dimensional compactum $Z$ with
  $\dim_{\z_p} Z \leq n$,
$\dim_{\z_{(q)}} Z \leq n$     and a cell-like  map $r : Z \lo X$ from $Z$  onto $X$.
\end{theorem}

This paper  goes along the line of investigation represented   by Theorems   \ref{dr1.5}, \ref{dr2}
and \ref{ko1}.   These theorems can be regarded as  particular cases of
 the following general problem:  Let $X$ be a compactum  with $\dim_\z X \leq n$.
 Do there exist an $n$-dimensional  compactum $Z $
 and a cell-like map from $Z$ onto $X$
 such that $\dim_G Z \leq \dim_G X$ for every abelian group $G$?
 The goal of this paper is to answer this problem affirmatively in cohomological
  dimensions
 larger than $1$.  Namely we will prove the following theorem.
 \begin{theorem}
 \label{t1}
 Let $X$ be a compactum  with $\dim_\z X \leq n \geq 2$.  Then there exist
 a compactum $Z$ with $\dim Z \leq n$  and a cell-like map $r : Z \lo X$ from
 $Z$ onto $X$ such that for every integer $k \geq 2$ and every group $G$  such that
 $\dim_G X \leq k$ we have $\dim_G Z \leq k$.
 \end{theorem}

 Theorem \ref{t1}  can be reformulated in terms of extensional dimension \cite{drdyd0, drdyd1}.
 The extensional dimension of $X$  is said not to exceed       a CW-complex  $K$,
  written $\edim X \leq K$,  if for every closed subset $A$
 of $X$ and every map $f : A \lo K$ there is an extension of $f$ over $X$.
 It is well-known that $\dim X \leq n$ is equivalent to $\edim X \leq \s^n$
 and $\dim_G X \leq  n$ is equivalent to $\edim X \leq K(G,n)$ where
 $K(G,n)$ is an Eilenberg-Mac Lane complex of type $(G,n)$.
 The following theorem  shows a close connection between  cohomological
 and extensional dimensions.
\begin{theorem}
 \label{dr3} {\rm \cite{dr3}}\qua
 Let $X$ be a
 compactum and let $K$ be a simply connected CW-complex. Consider
 the following conditions:

{\rm(1)}\qua $\edim X \leq  K$;

{\rm(2)}\qua $\dim_{H_i(K)} X \leq i$ for every $i > 1$;

{\rm(3)}\qua $\dim_{\pi_i(K)} X \leq i$ for every $i >1$.

Then   (2) and (3) are equivalent and (1) implies both (2) and (3).
If $X$ is  finite dimensional then all the conditions are equivalent.
\end{theorem}

Theorems \ref{t1} and \ref{dr3} imply the following:
\begin{theorem}
 \label{t2}
 Let $X$ be a compactum  with $\dim_\z X \leq n \geq 2$.  Then there exist
 a compactum $Z$ with $\dim Z \leq n$  and a cell-like map $r : Z \lo X$ from
 $Z$ onto $X$ such that for every  simply connected CW-complex $K$
   such that
 $\edim  X \leq K$ we have $\edim Z \leq K$.
 \end{theorem}
{\bf Proof}\qua Let $Z$ and $r : Z \lo X$ be as in Theorem \ref{t1}.
 Let a simply connected CW-complex $K$  be such that
 $\edim X \leq K$. Then by Theorem \ref{dr3},  $\dim_{H_i(K)} X \leq i$ for every $i > 1$
 and hence by Theorem  \ref{t1}, $\dim_{H_i(K)} Z \leq i$ for every $i > 1$. Then
 since $Z$ is finite dimensional  it follows from   Theorem \ref{dr3} that
 $\edim Z \leq K$.
\endproof

 Note that   the restriction $k \geq2$ in Theorem \ref{t1} cannot be omitted.
  Indeed, take an infinite dimensional compactum $X$ with $\dim_\q X=1$ and $\dim_\z X =2$
 (such a compactum was constructed by Dydak and Walsh \cite{d-w})
  and let $r : Z \lo X$ be a cell-like map of a $2$-dimensional compactum $Z$ onto $X$.
  Then  $\dim_\q Z =2 $  since otherwise  by
  a result of Daverman  \cite{da1} we would have      $\dim  X \leq 2$.
  This observation  also shows that Theorem \ref{t2} does not hold for
  non-simply connected complexes $K$.

  The author would like to thank the referee for the careful
  reading of the paper and valuable suggestions.

 \end{section}
 \begin{section}{Preliminaries}
     A  map between CW-complexes  is said to be  combinatorial if  the preimage of
   every subcomplex
   of the range is a subcomplex of the domain.

  Let $M$ be  a simplicial complex and let  $M^{[k]}$        be
    the $k$-skeleton of $M$ (=the union of all simplexes of $M$ of $\dim \leq k$).
By
 a resolution $EW(M,k)$   of $M$   we mean a CW-complex $EW(M,k)$ and
 a combinatorial map
 $\omega : EW(M,k) \lo  M$ such that $\omega$ is 1-to-1 over $M^{[k]}$.
 Let $f : N \lo K$  be a map of a subcomplex $N$ of $M$ into a CW-complex $K$.
The resolution is said to be suitable for $f$   if
 the map  $f \circ\omega|_{\omega^{-1}(N)}$ extends
   to a map  $ f': EW(M,k) \lo K$.  We will call $f'$  a  resolving map
   for $f$. The resolution is said to be
   suitable  for  a compactum $X$
if for   every simplex $\Delta$ of $M$,
 $\edim X  \leq \omega^{-1}(\Delta)$.
     Note that if $\omega: EW(M,k) \lo M$ is a resolution suitable
 for $X$   then  for every map  $\phi :  X \lo  M$  there is  a map   $\psi : X \lo EW(M,k)$
 such that  for every simplex $\Delta$ of $M$ ,
  $(\omega \circ \psi)(\phi^{-1}(\Delta)) \subset \Delta$.
 We will call $\psi$ a combinatorial lifting of $\phi$.

Let $M$   be a finite simplicial complex.
 Let       $f : N  \lo K$ be a cellular  map  from a subcomplex $N$
of $M$  to a CW-complex    $K$ such that
$M^{[k]}\subset N$.
 Following \cite{l0,l1}   we will  construct  a
 resolution   $\omega: EW(M,k) \lo M$  of
  $M$ which is
 suitable   for   $f $.
   In the sequel we will refer
 to this resolution as  the standard resolution for  $f$.
   We will associate     with the standard   resolution
 a  cellular resolving map $f' :  EW(M,k) \lo K$ which will be
 called the standard resolving map.
  The standard  resolution is constructed  by induction
 on $n=\dim (M\setminus N)$.

 For $M=N$ set   $EW (M,k)=  M$  and let    $\omega: EW(M,k) \lo M$
 be the identity map   with    the standard resolving map $f'=f$.
 Let $n > k$. Denote $M'=N\cup M^{[n-1]}$
 and   assume that     $ \omega': EW(M',k) \lo M'$  is  the standard
 resolution     of       $M' $
  for   $f$
  with the standard  resolving map  $f'  :  EW(M',k)  \lo K$.
     The standard resolution      $\omega: EW(M,k) \lo M$ is constructed
 as follows.

 The  CW-complex $EW(M,k)$  is obtained  from  $EW(M',k)$  by attaching
   the mapping cylinder of  $f'|_{{\omega'}^{-1}({\partial \Delta})}$
 to      ${\omega'}^{-1}({\partial \Delta})$  for every  $n$-simplex $\Delta$ of $M$
 which is not contained in $M'$.
 Let $\omega : EW(M,k) \lo M$ be the projection
which extends $\omega'$ by sending
  each mapping cylinder
 to  the corresponding $n$-simplex $\Delta$  such that
 the $K$-part of the cylinder
 is sent to the barycenter  of     $\Delta$ and   each interval connecting
 a point of ${\omega'}^{-1}({\partial \Delta})$
 with the corresponding point of the  $K$-part of the cylinder is
 sent linearly to
 the interval connecting the corresponding point
 of ${\partial \Delta}$ with the barycenter of  $\Delta$.
 We can naturally  define  the extension of
 $f'|_{ {\omega'}^{-1}({\partial \Delta})}$ over its mapping cylinder
 by sending   each interval of the cylinder to the corresponding point of $K$.
 Thus        we define the  standard resolving map
 which extends $f'$
      over $EW(M,k)$.     The CW-structure of
   $EW(M,k)$  is induced  by the CW-structure of   $EW(M',k)$ and
           the natural CW-structures of
   the     mapping cylinders in  $EW(M,k)$.
  Then with respect to this CW-structure
    the  standard resolving map is cellular
   and $\omega$ is combinatorial.

From the construction of the standard resolution it follows that
 for each simplex $\Delta$ of $M$,
    $\omega^{-1}(\Delta)$ is either contractible or
  homotopy equivalent to $K$ and   for every $x\in M$,
  $\omega^{-1}(x)$ is either a singleton or  homeomorphic to  $K$.
  It is  easy to check  that if $M$ and $K$ are $(k-1)$-connected
  then so is $EW(M,k)$.  Also note
that for every subcomplex $T$ of $M$,
  $\omega|_{\omega^{-1}(T)} :  EW(T,k)=  {\omega^{-1}}(T) \lo T$
  is the   standard resolution      of  $T$
  for $f|_{N \cap T}$.

 All  groups are assumed to be
   abelian  and functions between groups are homomorphisms.
   $\cal P$ stands for the set of  primes.  For a non-empty subset $\cal A$ of $ \cal P$ let
$S({\cal A})=\{ p_{1}^{n_1}p_{2}^{n_2}...p_{k}^{n_k} : p_i \in {\cal A}, n_i \geq 0\}$
 be  the set of  positive  integers with prime factors from $\cal A$ and for the empty set define
 $S(\emptyset)=\{ 1 \}$.
Let   $G$ be a group and  $g \in G$.
We say that $g$  is $\cal A$-torsion  if there is $n\in S(\cal A)$
such that $ng=0$  and  $g$  is $\cal A$-divisible if for every $n \in S(\cal A)$ there is $h\in G$
such that $nh=g$.
$\tor_{\cal A} G$ is the subgroup of the $\cal A$-torsion elements
of $G$.   $G$ is  $\cal A$-torsion  if $G=\tor_{\cal A}   G$, $G$ is $\cal A$-torsion free
if $\tor_{\cal A} G= 0$
and $G$ is $\cal A$-divisible if every element of $G$ is $\cal A$-divisible.

$G$ is $\cal A$-local if $G$ is $({\cal P} \sm {\cal A})$-divisible
and       $({\cal P} \sm {\cal A})$-torsion free.    The  $\cal A$-localization of  $G$ is
the homomorphism  $G \lo G \otimes \z_{(\cal A)}$  defined by $g \lo g \otimes 1$ where
$\z_{(\cal A)}= \{ n/m : n \in \z, m\in S({\cal P}\sm {\cal A})\}$.     $G$ is  $\cal A$-local
if and only if the $\cal A$-localization of  $G$  is an isomorphism.
A map between
two simply connected CW-complexes is  an $\cal A$-localization  if the induced homomorphisms of
the homotopy  and (reduced integral) homology groups
are $\cal    A$-localizations. 

Let $G$ be a group,  let $\alpha : L \lo M$ be a surjective  combinatorial map of a CW-complex
$L$ and a finite simplicial complex $M$ and let $n$ be a positive integer
 such that  $\Tilde{H}_i(\alpha^{-1}(\Delta) ; G)=0$
for every $i < n$ and  every  simplex $\Delta $ of $M$.     One can show by induction on
the number of simplexes of $M$  using  the Mayer-Vietoris sequence and the Five Lemma
that $\alpha_* :   \Tilde{H}_i (L; G) \lo    \Tilde{H}_i (M; G)$ is an isomorphism
for $i < n$.   We will refer to this fact as the combinatorial Vietoris-Begle theorem.

  \begin{proposition}
  \label{p1}
  Let $m \geq k+2$, $k \geq 2$ and let $M$ be an $(m-1)$-connected finite simplicial complex.
  Let $\omega : EW(M,k)\lo M$ be the standard resolution  for
 a cellular map  $f : N \lo K(G,k)$  from   a subcomplex $N$  of $M$ containing $M^{[k]}$.
  Then $EW(M,k)$ is $(k-1)$-connected and      for every $1\leq i \leq m-2$,
  $\pi_i (  EW(M,k))$ is

  {\rm(i)}\qua  $p$-torsion if $G=\z_p $;

  {\rm(ii)}\qua   $p$-torsion  and $\pi_k (  EW(M,k))$    is $p$-divisible
     if     $G=\z_{p^\infty} $;

  {\rm(iii)}\qua  $p$-local     if  $G=\z_{(p)} $   and $\emptyset$-local if $G=\q$.
  \end{proposition}
  {\bf Proof}\qua Since $M$ and $K(G,k)$ are $(k-1)$-connected then
  so is $EW(M,k)$.        Recall that
  $\omega$ is a surjective combinatorial  map and   for every simplex $ \Delta$ of  $ M$,
  $\omega^{-1}(\Delta)$ is either contractible  or homotopy equivalent to  $K(G,k)$.

\medskip
  {\bf (i)}\qua  By the generalized Hurewicz theorem the groups     $H_i ( K(\z_p, k))$, $i\geq 1$
  are $p$-torsion.  Then    $H_i ( K(\z_p, k))$, $i\geq 1$ is  $p$-local and
  $H_i( K(\z_p, k); \q)=H_i ( K(\z_p, k))\otimes\q =0$, $i \geq 1$.
  Let $q \in \cal  P$  and $q \neq p$.
  The $p$-locality of   $H_i ( K(\z_p, k))$, $i\geq 1$    implies that
          $H_i ( K(\z_p , k); \z_q)=0$, $i\geq 1$.
Then,  since $M$ is $(m-1)$-connected,    by  the combinatorial Vietoris-Begle   theorem
 we get  that    $   H_i(EW(M,k); \z_q)=0$  and  $H_i(EW(M,k); \q)=0$,
 $1\leq i \leq m-1$.
From the universal coefficient theorem it follows that  the last conditions imply
that     $H_i(EW(M,k))\otimes \q=0$
   for $1\leq i \leq m-1$     and    $H_i(EW(M,k)) * \z_q=0$
  for $1\leq i \leq m-2$.
 Hence $H_i(EW(M,k))$ is torsion  and $q$-torsion free for $1\leq i \leq m-2$
 and every        $q \in \cal  P$, $q \neq p$.
 Therefore            $  H_i(EW(M,k))$,  $1\leq i \leq m-2$  is
 $p$-torsion   and by the generalized Hurewicz theorem
  $  \pi_i (EW(M,k))$,  $1\leq i \leq m-2$   is       $p$-torsion.   

\medskip
 {\bf  (ii)}\qua
 Note that the proof of {\bf (i)} applies not only for $G=\z_p$ but also for
    $G=\z_{p^\infty}$.  Therefore
    we  can conclude  that
      $  \pi_i (EW(M,k))$  is
     $p$-torsion   for  $1\leq  i\leq  m-2$.

    By the Hurewicz theorem $\pi_k(EW(M,k))=   H_k(EW(M,k))$.  To show
    that      $H_k(EW(M,k))$ is $p$-divisible
     first observe that  $  H_k(K(\z_{p^\infty},k))=\z_{p^\infty}$ and
    by the universal coefficient theorem
    $H_k( K(\z_{p^\infty},k);\z_p)=  \z_{p^\infty} \otimes  \z_p=0$.
    Then since  $M$ is $k$-connected   the combinatorial  Vietoris-Begle theorem implies that
    $H_k (EM(M,k);\z_p)=0$.
    Once again by the universal coefficient
    theorem   $H_k( EW(M,k))\otimes \z_p =0$ and therefore    $H_k( EW(M,k))$ is $p$-divisible. 

\medskip
    {\bf (iii)}\qua
    We  will prove  the case $G=\z_{(p)}$.
    The case $G=\q$  is similar  to    $G=\z_{(p)}$.
     The proof applies  well-known results of
    Rational Homotopy Theory \cite{ratio}. The $p$-locality of  $\pi_i(K(\z_{(p)}, k))$
    implies that $ H_i(   K(\z_{(p)}, k))$, $i \geq 1$  are $p$-local.
    Then by the reasoning based on the combinatorial  Vietoris-Begle and  universal
coefficient theorems  that we used in the proof of {\bf (i)} one can show that
  $H_i (EW(M,k))$, $  1\leq i \leq m-2$   are  $p$-local and
  $H_{m-1}(EW(M,k))$ is $q$-divisible for every prime $q\neq p$.
  Let a map $ \alpha : EW(M,k)\lo L$  be
     a $p$-localization of $EW(M,k)$.  Then $\alpha $ induces an isomorphism
     of     $H_i (EW(M,k))$ and    $H_i(L)$  for   $  1\leq i \leq m-2$ and
     an epimorphism of  $H_{m-1} (EW(M,k))$ and    $H_{m-1}(L)$. Hence
     by the Whitehead theorem the groups $\pi_i (EW(M,k))$ and $\pi_i (L)$  are
     isomorphic   for  $  1\leq i \leq m-2$.   Thus       $\pi_i (EW(M,k))$, $  1\leq i \leq m-2$
     is $p$-local.
\endproof

     The following proposition  is  an infinite dimensional version of   the implication
  (3)$\Rightarrow$(1) of Theorem \ref{dr3}.

 \begin{proposition}
 \label{p2}
 Let $K$ be a simply connected CW-complex such that  $K$ has only finitely many non-trivial
 homotopy groups.   Let $X$ be a compactum such that $\dim_{\pi_i (K)} X  \leq i$
 for  $i>1$.  Then $\edim X \leq K$.
 \end{proposition}
 {\bf  Proof}\qua Let $n$ be such that $\pi_i(K)=0$ for $i\geq n$ and
 let $\psi: A\lo K$ be a map from  a closed subset $A$ of $X$ into $K$.
 Represent $X$ as the inverse limit $X=\invlim K_j$ of finite  simplicial complexes    $K_j$
 with combinatorial bonding maps
 $h^{j}_{l} : K_{l} \lo K_{j}, l > j$  and the  projections  $p_j : X \lo K_j$ such that
 diam$p_j^{-1}(\Delta) \leq 1/j$ for every simplex  $\Delta $ of $K_j$.
 Let $j$ be so large  that there is a map $f :  N \lo K$  from a subcomplex $N$ of $K_j$
 such that $A \subset p_j^{-1}(N)$ and  $f \circ p_j |_A$ is homotopic to $\psi$.
 Then,   since $ {\check H}^{i+1}(X,  p_j^{-1}(N)  ; \pi_i (K))=0$   for every $i$,
  by Obstruction Theory there is       a sufficiently large $l  > j$ such that
    $f \circ  h^{j}_{l}  |_{(h^j_{l})^{-1}(N)}$
   extends over   the $n$-skeleton of    $K_{l}$ and since
 $\pi_{i}(K)=0$ for $i\geq n$   the map    $f \circ  h^j_{l}  |_{(h^j_{l})^{-1}(N)}$      will also extend
 over $K_{l}$ to $f' : K_{l} \lo K$.  Then $f' \circ p_{l}|_A$ is homotopic to $\psi$ and hence
 $\psi$ extends over $X$.  Thus  $\edim X \leq K$.
\endproof

Let $K'$ be a simplicial complex.  We say that   maps
  $h : K \lo K'$, $g : L \lo L'$,  $\alpha : L \lo K$ and $\alpha' : L' \lo K'$
  \[
    \begin{CD}
           L @>\alpha>>   K\\
           @V g VV        @V h VV \\
             L' @>\alpha'>>   K'
    \end{CD}
\]
 combinatorially commute   if   for every
 simplex  $\Delta$ of $K'$ we have that
 $(\alpha' \circ g)(( h\circ \alpha)^{-1}(\Delta)) \subset \Delta$.
(The direction in which we want the maps $ h, g, \alpha$ and $\alpha'$
 to   combinatorially
 commute    is indicated by the first map  in the list. Thus saying
 that $\alpha',  h, g $ and $\alpha$ combinatorially commute we would mean that
   $(h\circ \alpha)(( \alpha'\circ  g)^{-1}(\Delta)) \subset \Delta$
   for every simplex $\Delta $ of $K'$.)
A map $h' : K \lo L'$ is said to be
 a combinatorial lifting of $h$  to $L'$  if
 for every  simplex  $\Delta$ of $K'$   we have that
 $(\alpha' \circ h')( h^{-1}(\Delta)) \subset \Delta$.

 For   a simplicial complex $K$ and $a \in K$,   $st(a)$  denotes the union of all the simplexes
 of   $K$ containing  $a$.
 The following proposition    whose
 proof  is left to the reader      is
 a collection of   simple combinatorial  properties of    maps.

 \begin{proposition}
 \label{p3}${}$

   {\rm(i)}\qua     Let a compactum $X$  be  represented as the inverse limit
 $X  ={\rm  \invlim }K_i$ of  finite simplicial complexes $K_i$
 with  bonding maps   $h_{j}^i  :  K_{j} \lo K_i$.  Fix $i$ and   let
 $\omega : EW(K_i,k) \lo K_i$ be a resolution of $K_i$ which is suitable for $X$.
 Then there is
  a sufficiently large $j$  such that $h_j^i$ admits a combinatorial lifting to $EW(K_i,k)$.

  {\rm(ii)}\qua        Let  $h : K \lo K'$, $h' : K \lo L'$ and $\alpha' : L' \lo K'$ be maps of
  a simplicial complex $K'$   and  CW-complexes $K$ and $L'$  such that
  $h$ and $\alpha'$ are combinatorial and $h'$ is a combinatorial lifting of $h$.
 Then  there is a cellular approximation of $h'$ which is  also a combinatorial lifting of $h$.

 {\rm(iii)}\qua     Let $K$ and $K'$ be simplicial complexes,
 let    maps    $h : K \lo K'$, $g : L \lo L'$,  $\alpha : L \lo K$ and $\alpha' : L' \lo K'$
 combinatorially commute and let $h$ be  combinatorial.  Then

 $ g( \alpha^{-1}(st( x)) ) \subset
   \alpha'{}^{-1}(st( h (x)) )$ and
   $h(st(\alpha (z)))  \subset        st((\alpha' \circ g)(z))$ 

   for every   $x \in K$ and  $z\in L$.

 \end{proposition}

  We end this section with recalling  basic facts of Bockstein Theory.
  The Bockstein basis  is
 the following collection  of groups $\sigma = \{\q,  \z_p,   \z_{p^\infty}, \z_{(p)}    : p \in {\cal P} \}$.
 For an abelian group $G$  the Bockstein basis
 $\sigma(G)$  of $G$
    is  a subcollection  of $\sigma$ defined as follows:   

    $\z_{(p)} \in \sigma(G)$  if   $G/\tor G$ is not divisible by $p$;

    $\z_p \in    \sigma(G)$ if $\tor_p G $ is not divisible by $p$;

    $\z_{p^\infty} \in   \sigma(G)$  if     $\tor_p G \neq 0$  and  $\tor_p G $  is divisible by $p$;

    $ \q \in        \sigma(G)$ if      $G/\tor G\neq 0$  and     $G/\tor G$ is divisible by every $p \in \cal P$.

  Let $X$ be a compactum.   The Bockstein theorem says that
 
 $  \dim_G X = \sup\{\dim_E X : E \in       \sigma(G) \}$.

  The Bockstein inequalities relate
  the cohomological dimensions for different groups of Bockstein basis.
  We will use   the  following inequalities:  

  $\dim_{\z_{p^\infty}} X \leq
  \dim_{\z_p} X \leq     \dim_{\z_{p^\infty}} X +1$;

  $\dim_{ \z_p} X \leq     \dim_{ \z_{(p)} }X$ and     $\dim_{ \q} X \leq     \dim_{ \z_{(p)}} X$.
 
  Finally  recall  that $\dim_G X \leq \dim_\z X$ for every abelian group $G$.

  \end{section}

  \begin{section}{Proof of Theorem \ref{t1}}
  Let $m=n+2$.
    Represent $X$ as   the inverse limit $X  = \invlim (K_i,h_i)$ of finite simplicial complexes $K_i$
 with combinatorial bonding maps   $h_{i+1}  :  K_{i+1} \lo K_i$   and the projections
 $p_i : X \lo K_i$ such that for every simplex $\Delta$ of $K_i$,  diam$(p_i^{-1}(\Delta)) \leq 1/i$.
  We will construct  by induction finite  simplicial complexes
  $L_i$   and maps      $g_{i+1}: L_{i+1} \lo L_i$,
 $\alpha_i : L_i \lo K_i$   such that

 (a)\qua  $L_i= K_i^{[m]}$   and $\alpha_i    : L_i \lo K_i$ is the inclusion.  The simplicial structure
 of $L_1$ is      induced from     $ K_1^{[m]}$         and
  the simplicial structure
 of $L_i$, $i >1$ is
   defined as   a  sufficiently  small barycentric subdivision of   $ K_i^{[m]}$.
   We will refer to this simplicial structure while constructing standard resolutions
   of $L_i$.
    It is clear that   $\alpha_i$ is
 always a combinatorial map;

 (b)\qua  the maps $h_{i+1}$, $g_{i+1}$,  $\alpha_{i+1}$  and $\alpha_i$   combinatorially commute.
 Recall that   this  means that
 for every simplex  $\Delta $ of $ K_{i}$,
 $(\alpha_i  \circ g_{i+1})((h_{i+1} \circ \alpha_{i+1})^{-1}(\Delta)) \subset \Delta$.

  We will construct $L_i$ in such a way that   $Z=\invlim(L_i,g_i)$ will be of $\dim  \leq n$
  and $Z$  will admit a cell-like  map  onto $X$ satisfying the conclusions of the theorem.
 Assume that the construction
 is completed for $i$.  We  proceed to $i+1$ as follows.

   Let $E \in \sigma$ be such that
 $\dim_E X \leq k$, $2 \leq k  \leq n$ and let $f  :  N  \lo K(E,k)$ be a  cellular map  from
 a subcomplex $N$ of $L_i$, $  L_i^{[k]}  \subset  N$.
Let  $\omega_L: EW(L_i, k) \lo L_i^{}$
be the standard resolution  of $L_i$ for $f$.   We are  going to construct
from  $\omega_L: EW(L_i, k) \lo L_i^{}$ a resolution
 $\omega : EW(K_i, k) \lo K_i$
 of $K_i$ suitable for $X$.
   If $\dim K_i \leq  k$ set
 $\omega=\alpha_i \circ \omega_L      :  EW(K_i,k)=EW(L_i, k) \lo K_i$.

If $\dim K_i  >  k$ set
$\omega_{k}=\alpha_i \circ \omega_L :     EW_{k}(K_i,k)=EW(L_i, k) \lo  K_i$
  and we will construct by induction  resolutions
 $\omega_j : EW_j(K_i, k) \lo K_i$, $k+1\leq j\leq \dim K_i$
  such that   $EW_{j}(K_i, k)$ is a subcomplex of
   $EW_{j+1}(K_i, k)$ and $\omega_{j+1}$ extends $\omega_j$  for every $k\leq j  < \dim K_i$.
    The construction  is carried out as follows.

    Assume that $\omega_j :   EW_j(K_i, k) \lo K_i$, $k\leq j <\dim K_i$
    is constructed.
    For every simplex $\Delta$ of $K_i$ of $\dim=j+1 $   consider  the subcomplex
    ${\omega_j}^{-1}(\Delta)$ of  $EW_j(K_i, k)$.  Enlarge
    ${\omega_j}^{-1}( \Delta)$    by
  attaching  cells  of $\dim \geq m+1  $
   in order     to get a subcomplex with trivial homotopy groups in
 $\dim \geq  m$.
 Let   $EW_{j+1}(K_i, k)$ be $EW_j(K_i, k) $ with all the cells  attached   for
 all  $(j+1)$-dimensional simplexes  $\Delta$ of $K_i$ and let
  $\omega_{j+1}:      EW_{j+1} (K_i, k)\lo K_i$ be an extension of
  $\omega_j$
      sending the interior points of the attached cells to the interior of       the corresponding
    $\Delta$.

      Finally denote  $EW(K_i, k)=EW_j(K_i, k)$ and $\omega=\omega_j : EW_j(K_i, k) \lo K_i$
    for $j=\dim K_i$.       Note that since we attach cells only of $\dim >m$,
    the $m$-skeleton of   $EW(K_i, k)$     coincides   with  the $m$-skeleton of $EW(L_i,k)$.

        Let us show that $EW(K_i, k)$ is suitable for $X$.
        Fix a simplex $\Delta$ of $K_i$ and denote $M= \alpha_i^{-1}(\Delta)$.
       First note that  $M$  is $(m-1)$-connected,
       $\omega^{-1}(\Delta )$  is $(k-1)$-connected,
          $\pi_j(\omega^{-1}(\Delta ))=0$  for $j \geq  m$
    and $\pi_j(\omega^{-1}(\Delta ))=\pi_j (\omega_L^{-1}(M))$ for $k \leq j \leq n$.
   Also note that  since $\dim_\z X \leq n$,
         $\dim_{ \pi_j ( \omega^{-1}(\Delta ))} X \leq  \dim_{\z} X \leq n$    for   $n \leq  j$.
In order to show that $\edim X \leq  \omega^{-1}(\Delta)$
    consider   separately the following cases. 

    {\bf Case 1}\qua $E=\z_p$.  By (i) of  Proposition \ref{p1},
        $\pi_j(\omega^{-1}(\Delta ))=\pi_j (\omega_L^{-1}(M))$, $k \leq j \leq n$
        is $p$-torsion. Hence
        by Bockstein Theory
        $\dim_{ \pi_j ( \omega^{-1}(\Delta ))} X \leq  \dim_{\z_p} X \leq k$  for  $k\leq j  \leq n$.
        Therefore  by Proposition \ref{p2}, $\edim X \leq     \omega^{-1}(\Delta )$.

    {\bf Case 2}\qua $E=\z_{p^\infty}$.     By (ii) of  Proposition \ref{p1},
        $\pi_j(\omega^{-1}(\Delta ))=\pi_j (\omega_L^{-1}(M))$, $k \leq j \leq n$ is $p$-torsion  and
               $\pi_k(\omega^{-1}(\Delta ))=\pi_k (\omega_L^{-1}(M))$ is $p$-divisible.   Hence
        by the Bockstein       theorem and inequalities
          $\dim_{ \pi_k ( \omega^{-1}(\Delta ))} X \leq  \dim _{\z_{p^\infty}} X \leq k$
          and
$\dim_{ \pi_j ( \omega^{-1}(\Delta ))} X \leq  \dim_{\z_p} X \leq \dim _{\z_{p^\infty}} X +1 \leq k+1$
          for $k+1 \leq j \leq n $.
                  Therefore by
         Proposition \ref{p2}, $\edim X \leq     \omega^{-1}(\Delta )$.

 {\bf Case 3}\qua $E=\z_{(p)}$.  By (iii) of  Proposition \ref{p1},
        $\pi_j(\omega^{-1}(\Delta ))=\pi_j (\omega_L^{-1}(M))$, $k \leq j \leq n$
        is $p$-local.
Then $\sigma(   \pi_j(\omega^{-1}(\Delta )))$ may possibly contain only the groups
$\z_p$, $\z_{p^\infty}$,  $\z_{(p)}$ and $\q$.
        Hence
        by the Bockstein  theorem and inequalities
        $\dim_{ \pi_j ( \omega^{-1}(\Delta ))} X \leq  \dim_{\z_{(p)}} X \leq k$  for every $k\leq j \leq  n$.
 Therefore by Proposition \ref{p2}, $\edim X \leq     \omega^{-1}(\Delta )$.

        {\bf  Case 4}\qua  $E=\q$.       By (iii) of  Proposition \ref{p1},
        $\pi_j(\omega^{-1}(\Delta ))=\pi_j (\omega_L^{-1}(M))$, $k \leq j \leq n$
        is $\emptyset$-local.    Then $\sigma(   \pi_j(\omega^{-1}(\Delta )))$ may possibly
        contain only $\q$ and hence
          $\dim_{ \pi_j ( \omega^{-1}(\Delta ))} X \leq  \dim_{\q} X \leq k$  for every $k\leq j \leq n$.
 Therefore by Proposition \ref{p2}, $\edim X \leq     \omega^{-1}(\Delta )$.    

          Thus we have shown that $EW(K_i, k)$ is suitable for $X$.
          Replacing  $K_{i+1}$ by
         $K_j$ with a sufficiently large $j$ we may   assume
    by (i)  of Proposition  \ref{p3}       that there is a combinatorial lifting
         of $h_{i+1}$ to $h'_{i+1} : K_{i+1} \lo  EW(K_i,k)$.
        By (ii) of  Proposition  \ref{p3} we  replace $h'_{i+1}$ by its  cellular approximation
         preserving
  the property of   $h'_{i+1}$ of being a      combinatorial lifting of $h_{i+1}$.

Then $h'_{i+1}$  sends the $m$-skeleton of $K_{i+1}$ to
the $m$-skeleton of $EW(K_i,k)$.  Recall that the  $m$-skeleton of $EW(K_i,k)$ is
contained in  $EW(L_i,k)$ and hence
 one can define
 $g_{i+1}=\omega_L \circ h'_{i+1}|_{K_{i+1}^{[m]}} : L_{i+1}= K_{i+1}^{[m]}\lo L_i$.
 Finally define a simplicial structure on $L_{i+1}$ to be a sufficiently small barycentric
 subdivision of      $K_{i+1}^{[m]}$  such that   

  (c)\qua
 diam$g_{i+1}^j(\Delta) \leq 1/i$ for every simplex
 $\Delta$ in  $L_{i+1}$  and $j \leq i$   
     where
 $g^j_i=g_{j+1} \circ g_{j+2} \circ ...\circ g_i: L_i \lo L_j$.   

  It is easy to check that the properties (a) and (b)
 are satisfied.

        Denote $Z=\invlim(L_i,g_i)$ and let $r_i : Z \lo L_i$ be the projections.
            For constructing
           $L_{i+1}$  we used an arbitrary map $f : N  \lo K(E,k), E\in \sigma$ such that
           $\dim_E X \leq k$,  $ 2 \leq k \leq n$  and  $N$ is a subcomplex
 of $L_i$ containing $L_i^{[k]}$.
            Let us show that  choosing $E \in \sigma$ and $f$  in an appropriate way
           for each $i$ we can achieve
that $\dim_G Z \leq k$  for every integer $k$ and
 group $G  $  such that      $\dim_G  X \leq k \geq 2$.

  Let $E\in \sigma$ be such that $\dim_E X \leq k $,  $ 2 \leq k \leq n$
  and let $\psi : F \lo K(E,k)$ be a map of a closed subset
 $F$ of $L_j$.
  Then by (c) for a sufficiently large $i>j$ the map
 $\psi \circ g_i^j|_{(g_i^j)^{-1}(F)}$ extends over a subcomplex $N$ of $L_i$
to  a map $f : N \lo  K(E,k)$.  Extending
    $f$ over  $L_i^{[k]}$ we may assume that $   L_i^{[k]} \subset N$.
       Replacing      $f$ by  its cellular approximation we
    also assume that  $f$ is cellular.          Now suppose that we use this   map $f$
   for  constructing  $L_{i+1}$.

          Since $g_{i+1}$ factors
           through $  EW(L_i, k)$, the   map $f \circ g_{i+1}|_{g^{-1}_{i+1}(N)}:
 {g^{-1}_{i+1}(N)} \lo K(E, k)$ extends to a map $f' : L_{i+1}\lo K(E, k)$.
Then $f'|_{(g^{j}_{i+1})^{-1}(F)}$ is homotopic to
  $\psi \circ g_{i+1}^j |_{(g^{j}_{i+1})^{-1}(F)}:
 {(g^{j}_{i+1})^{-1}(F)} \lo K(E, k)$   and therefore
   $\psi \circ g_{i+1}^j |_{(g^{j}_{i+1})^{-1}(F)}$      extends over $L_{i+1}$.
Now since we need to solve only countably many extension problems
for every $L_j$ with respect to $K(E,k)$ for every $E \in \sigma$
such that $\dim_E X \leq k$, $ 2\leq k \leq n$
we can choose
  for each $i$ a map $f: N  \lo  K(E,k)$  from a subcomplex $N$ of $L_i$   in the way
  described above to achieve that  $\dim_E Z \leq k$ for every $E \in \sigma$ such that
   $  \dim_E X \leq k ,  2\leq k   \leq n$ (an algorithm how to
   assign extension problems to various indices $i$ can be found
   in \cite{dr2} and \cite{dr-rep}).
   Then        by Bockstein  Theory    $\dim_\z  Z \leq n$ and
    $\dim_G Z \leq k$ for every  $G$ such that
   $  \dim_G X \leq k   \geq 2$.
   Since $Z$ is finite dimensional and  $\dim_\z  Z \leq n$  we get that   $\dim Z \leq n$.

   By (iii) of Proposition \ref{p3},        the properties (a) and   (b) imply
    that for every   $x \in X$ and  $z\in Z$
   the following  holds:

 (d1)
  $ g_{i+1}( \alpha_{i+1}^{-1}(st( p_{i+1} (x) ) )) \subset
   \alpha_{i}^{-1}(st( p_i (x)) )$  and

   (d2)
   $h_{i+1}(st((\alpha_{i+1}\circ r_{i+1})(z)))  \subset        st((\alpha_{i}\circ r_{i})(z))$.

   Define a  map $r : Z \lo X$ by $r(z)=\cap \{ p_i^{-1}(    st((\alpha_{i}\circ r_{i})(z) ) ): i=1,2,... \}$.
  Then   (d2) implies
   that $r$ is indeed well-defined and continuous.

  The properties (d1) and (d2) also imply  that  for every $x \in X$  

 $r^{-1}(x)=\invlim ( \alpha_i ^{-1}(st(p_{i} (x))), g_i |_{\alpha_i ^{-1}(st(p_{i} (x)))})$

 where  the map  $ g_i |_{...} $ is considered as a map
      to $\alpha_{i-1} ^{-1}(st(p_{i-1} (x)))$. 

Since $r^{-1}(x) $ is not empty for every $x \in X$,
 $r$ is a map onto and let us show that $r^{-1}(x)$ is
cell-like.        Let $ \phi :  r^{-1}(x) \lo K$ be a map to a CW-complex $K$. Then since
  $r^{-1}(x)=\invlim ( \alpha_i ^{-1}(st(p_{i} (x))), g_i |_{...})$  there is   a sufficiently large $i$
   such that
  the map $\phi$ can be factored up to homotopy through  the map
  $\gamma=r_i |_{r^{-1}(x)}:   r^{-1}(x) \lo T=\alpha_i ^{-1}(st(p_{i} (x)))$, that is
  there is  a map
  $\beta :   T  \lo K$   such that
   that  $\phi$ is homotopic
  to $\beta \circ \gamma$.
  Note that $st(p_{i} (x))$ is contractible and
  $T$ is homeomorphic to   the  $m$-skeleton of  $st(p_{i} (x))$.  Hence
  $T$ is    $(m-1)$-connected and since $  r^{-1}(x)$ is  of $\dim \leq n =m-2$,
  the map $\gamma$ is  null-homotopic.  Then  $\phi$ is also null-homotopic and hence
   $r$ is a cell-like map.  The theorem is proved.
\endproof

  \end{section}

\Addresses\recd
\end{document}

%% file: agtout.tex
\def\ifplaintex{\expandafter\ifx\csname documentclass\endcsname\relax}

\def\gtp{{\mathsurround=0pt\it $\cal G\mskip-2mu$eometry \&\ 
$\cal T\!\!$opology $\cal P\!$ublications}}  

\def\recd{{\small Received:\qua\receiveddate\ifx\reviseddate\relax
\else\qquad Revised:\qua\reviseddate\fi\par}} 

\def\lognumber#1{\def\thelognumber{#1}}
\def\volumenumber#1{\def\thevolumenumber{#1}}
\def\volumeyear#1{\def\thevolumeyear{#1}}
\def\papernumber#1{\def\thepapernumber{#1}}
\def\pagenumbers#1#2{\def\startpage{#1}\def\finishpage{#2}}
\def\published#1{\def\publishdate{#1}}

\def\received#1{\def\receiveddate{#1}}
\def\revised#1{\def\reviseddate{#1}}
\def\accepted#1{\def\accepteddate{#1}}

\long\def\asciiabstract#1{\long\def\theasciiabstract{#1}}

\let\\\par\let\thelognumber\relax\let\thevolumenumber\relax
\let\thepapernumber\relax\let\thevolumeyear\relax\let\startpage\relax
\let\finishpage\relax\let\publishdate\relax\let\receiveddate\relax
\let\reviseddate\relax\let\accepteddate\relax\let\theasciititle\relax
\let\theasciiauthors\relax
\let\theasciiabstract\relax

\let\theasciiemail\relax

\ifplaintex
\font\logobig=cmssbx10 scaled 3836
\font\logomed=cmssbx10 scaled 2557
\else
\font\logobig=cmssbx10 scaled 4200
\font\logomed=cmssbx10 scaled 2800
\fi

\long\def\makeagttitle{   
\count0=\startpage
\agt\hfill      
\hbox to 45truept{\vbox to 0pt{\vglue -13truept{\logomed A\kern -.37em{\logobig 
T}\kern -.38em G}\vss}\hss}
\break
{\small Volume \thevolumenumber\ (\thevolumeyear)
\startpage--\finishpage\nl
Published: \publishdate}

\vglue .25truein

{\parskip=0pt\leftskip 0pt plus
1fil\def\\{\par\smallskip}{\Large\bf\thetitle}\par\medskip} \vglue
0.05truein

{\parskip=0pt\leftskip 0pt plus 1fil\def\\{\par}{\sc\theauthors}
\par\medskip}%
 
\vglue 0.03truein 

{\small\leftskip 25truept\rightskip 25truept{\bf Abstract}\stdspace\theabstract

{\bf AMS Classification}\stdspace\theprimaryclass
\ifx\thesecondaryclass\relax\else; \thesecondaryclass\fi\par
{\bf Keywords}\stdspace \thekeywords\par}\vglue 7truept

}   

\ifplaintex
\font\phead=cmsl9 scaled 950
\font\pnum=cmbx10 scaled 913
\font\pfoot=cmsl9 scaled 950
\headline{\vbox to 0pt{\vskip -4.5mm\line{\small\phead\ifnum
\count0=\startpage ISSN 1472-2739 (on-line) 1472-2747 (printed)
\hfill {\pnum\folio}\else\ifodd\count0\def\\{ }%
\ifx\theshorttitle\relax\thetitle\else\theshorttitle\fi\hfill{\pnum\folio}
\else\def\\{ and }{\pnum\folio}\hfill\ifx\theshortauthors\relax\theauthors
\else\theshortauthors\fi\fi\fi}\vss}}
\footline{\vbox to 0pt{\vglue 0mm\line{\small\pfoot\ifnum\count0=\startpage
\copyright\ \gtp\hfill\else
\agt, Volume \thevolumenumber\ (\thevolumeyear)\hfill\fi}\vss}}
\else
\font=cmsl9 scaled 1050
\font\lnum=cmbx10 
\font=cmsl9 scaled 1050
\makeatletter
\def\@oddhead{{\small\lhead\ifnum\count0=\startpage ISSN 1472-2739 
(on-line) 1472-2747 (printed)\hfill {\lnum\number\count0}\else\ifodd\count0
\def\\{ }\ifx\theshorttitle\relax \thetitle \else\theshorttitle\fi\hfill
{\lnum\number\count0}\else\def\\{ and }{\lnum\number\count0}
\hfill\ifx\theshortauthors\relax 
\theauthors\else\theshortauthors\fi\fi\fi}}\def\@evenhead{\@oddhead}
\def\@oddfoot{\small\lfoot\ifnum\count0=\startpage\copyright\ \gtp\hfill\else
\agt, Volume \thevolumenumber\ (\thevolumeyear)\hfill\fi}
\def\@evenfoot{\@oddfoot}
\makeatother
\fi
\let\maketitlepage\makeagttitle
\let\makeshorttitle\maketitlepage
\let\maketitle\maketitlepage

\newwrite\gtoutfile
\long\gdef\makeheadfile{  
{\def\\{, }\def\s{ }
\immediate\openout\gtoutfile head.xxx
\immediate\write\gtoutfile{Proxy-for: \ifx\theasciiauthors\relax
\theauthors\else\theasciiauthors\fi\s<\ifx\theasciiemail\relax\theemail\else\theasciiemail\fi>}
\immediate\write\gtoutfile{\noexpand\\}
\immediate\write\gtoutfile{Authors: \ifx\theasciiauthors\relax
\theauthors\else\theasciiauthors\fi}
{\def\\{ }\immediate\write\gtoutfile{Title: \ifx\theasciititle\relax
\thetitle\else\theasciititle\fi}}
\immediate\write\gtoutfile{Subj-class: GT or SG, GR etc}
\immediate\write\gtoutfile{MSC-class: \theprimaryclass\ifx\thesecondaryclass\relax\else, \thesecondaryclass\fi}
\immediate\write\gtoutfile{Journal-ref: Algebr. Geom. Topol. \thevolumenumber\s
(\thevolumeyear) \startpage-\finishpage}
\immediate\write\gtoutfile{Comments: Published by Algebraic and
Geometric Topology at}
\immediate\write\gtoutfile{\s\s\s  http://www.maths.warwick.ac.uk/agt/AGTVol\thevolumenumber/agt-\thevolumenumber-\thepapernumber.abs.html}
\immediate\write\gtoutfile{\noexpand\\}
\immediate\write\gtoutfile{}
\ifx\theasciiabstract\relax
\immediate\write\gtoutfile{\theabstract}\else
\immediate\write\gtoutfile{\theasciiabstract}\fi
\immediate\write\gtoutfile{}
\immediate\write\gtoutfile{\noexpand\\}
\immediate\write\gtoutfile{}
\immediate\closeout\gtoutfile}}  

\def\maketitlepage{\makeagttitle\makeheadfile}
\let\makeshorttitle\maketitlepage
\let\maketitle\maketitlepage

%% file: agt-3-46.bbl
\begin{thebibliography}{99}
 \bibitem{da1}
R. J.  Daverman, \emph{ Hereditarily aspherical compacta and
cell-like maps.} Topology Appl. 41 (1991),
   no. 3, 247--254.
\bibitem{dr0} A. N. Dranishnikov, \emph{ On a problem of P. S. Aleksandrov.} (Russian)
Mat. Sb. (N.S.) 135(177) (1988),
no. 4, 551--557, 560; translation in Math. USSR-Sb. 63 (1989), no. 2, 539--545.
 \bibitem{dr1}
A. N.   Dranishnikov, \emph{ On homological dimension modulo $p$.}
   (Russian) Mat. Sb. (N.S.) 132(174)
(1987), no. 3, 420--433, 446; translation in Math. USSR-Sb. 60 (1988), no. 2, 413--425.
\bibitem{dr1.5}
A. N.  Dranishnikov, \emph{ $K$-theory of Eilenberg-Mac Lane
spaces and cell-like mapping problem.} Trans. Amer. Math. Soc. 335
(1993), no. 1, 91--103.
\bibitem{dr3} A. N.   Dranishnikov, \emph{ Extension of mappings into CW-complexes.}
(Russian) Mat. Sb. 182 (1991), no. 9, 1300--1310;
translation in Math. USSR-Sb. 74 (1993), no. 1, 47--56
 \bibitem{dr2} A. Dranishnikov, \emph{ Cohomological dimension theory of compact metric spaces.}
  Topology Atlas Invited Contributions,   http://at.yorku.ca/topology/taic.htm

  \bibitem{drdyd0}  A.  Dranishnikov, J. Dydak,
 \emph{ Extension dimension and extension types.} Tr. Mat. Inst. Steklova 212
(1996), Otobrazh. i Razmer., 61--94; translation in Proc. Steklov Inst. Math. 1996, no. 1 (212), 55--88
   \bibitem{drdyd1}
 Alexander Dranishnikov, Jerzy Dydak, \emph{ Extension theory of separable metrizable spaces with
applications to dimension theory.} Trans. Amer. Math. Soc. 353
(2001), no. 1, 133--156.
\bibitem{dr-rep}
  Alexander N. Dranishnikov,   Du\v san Repov\v s,
  \emph{ On Alexandroff theorem for general abelian groups.}
   Topology Appl. 111 (2001), no. 3, 343--353.
\bibitem{ed1}   R. D. Edwards, \emph{A theorem  and a question  related to cohomological
dimension and cell-like maps.} Notices of  the AMS,  25(1978),
A-259.
\bibitem{ratio}     Yves Felix, Stephen Halperin, Jean-Claude Thomas,
 \emph{  Rational homotopy theory.}   Graduate texts in mathematics no. 205,
   New York,  Springer-Verlag, 2001.
 \bibitem{d-w}
  Jerzy Dydak, John J. Walsh, \emph{ Infinite-dimensional compacta having cohomological dimension two:
   an application of the Sullivan conjecture.} Topology 32 (1993), no. 1, 93--104.
  \bibitem{ko1}  Akira Koyama, Katsuya Yokoi,
  \emph{On Dranishnikov's cell-like resolution.} Geometric topology:
Dubrovnik 1998. Topology Appl. 113 (2001), no. 1-3, 87--106.
 \bibitem{l0}    Michael Levin, \emph{ Constructing compacta of different extensional dimensions.}
 Canad. Math. Bull. 44(2001), no. 1, 80--86.
 \bibitem{l1}
 Michael Levin,\emph{  Acyclic resolutions for arbitrary groups.}
 Israel J. Math. 135 (2003), 193--203.

\bibitem{w1}
 John J. Walsh,
 \emph{ Dimension, cohomological dimension, and cell-like mappings.}
  Shape theory and
geometric topology (Dubrovnik, 1981), pp. 105--118,
Lecture Notes in Math., 870, Springer, Berlin-New York, 1981.

   \end{thebibliography}
